\documentclass[12pt]{article}
\usepackage{amsmath,amssymb,theorem,comment}
{\theorembodyfont{\rmfamily} 
\newtheorem{Rem}{Remark}}
{\theorembodyfont{\rmfamily} 
\newtheorem{Ex}{Example}}
{\theorembodyfont{\rmfamily} 
}
{\theorembodyfont{\rmfamily} 
}
\newtheorem{Thm}{Theorem}
\newtheorem{Lem}{Lemma}
\newtheorem{Prop}{Proposition}

\newenvironment{proof}{\noindent{\em proof.}}{\hfill $\square$ \medskip}
\makeatletter
\def\address#1#2{\begingroup
\noindent\parbox[t]{7.8cm}{%
\small{\scshape\ignorespaces#1}\par\vskip1ex
\noindent\small{\itshape E-mail address}%
\/: #2\par\vskip4ex}\hfill%
\endgroup}%
\makeatother

\title{
A limit transition from the Heckman-Opdam hypergeometric 
functions to the Whittaker functions associated with root systems
}
\date{}
\author{Nobukazu Shimeno}

\begin{document}
\maketitle
\footnote{ 
2000 {\em Mathematics Subject Classification}.
Primary 22E30, 33C67; Secondary 43A90.
}

\begin{abstract}
We prove that  the radial part of the class one Whittaker function on 
a split semisimple Lie group can be obtained as an appropriate limit of the 
Heckman-Opdam hypergeometric function. 
\end{abstract}

\section*{Introduction}
Among quantum integrable systems associated with root 
systems, there are three classes where 
there are well behaved 
joint eigenfunctions  that are   
closely related with Lie theory. They are 
the trigonometric Calogero-Moser model, the rational CM model, 
and the Toda model. For the CM model, the eigenfunction is 
the Heckman-Opdam hypergeometric function and 
the Bessel function corresponding to  trigonometric and 
rational cases respectively (cf. \cite{heckman,opdam0,opdam}). 
Among other eigenfunctions, they are up to constant multiples 
unique globally defined analytic 
functions. 
For special parameter they are the radial part of the 
spherical functions on a Riemannian symmetric space of the 
non-compact type and the Euclidean type respectively (cf. \cite{helgason2}). 
For the non-periodic Toda model, the eigenfunction is  
the class one Whittaker function defined by 
the Jacquet integral  on a Riemannian symmetric space of the 
non-compact type (cf. \cite{jacquet,hashizume}). 
Among other eigenfunctions, it is 
 up to a constant multiple unique eigenfunction of moderate growth. 

On the other hand, there are two limit transitions between the Hamiltonians,  
one is 
from the trigonometric CM model to the rational CM model, and 
the other is from the trigonometric CM model to the Toda model 
(cf. \cite{deJeu, inozemtsev}). 
In rank one case, corresponding limit transitions are one 
from the Gauss hypergeometric function to the Bessel function, 
and the other is from the Gauss hypergeometric function to 
the Macdonald function. In general, 
a limit transition for eigenfunctions in the former case was 
established by  Ben~Sa\"{i}d-\O rsted~\cite{BO2} and de Jeu~\cite{deJeu}. 
In this paper we establish a limit transition in the latter case. Namely 
we prove that a limit of the Heckman-Opdam hypergeometric function 
is the radial part of the Whittaker function on a split semisimple Lie group 
(Theorem~\ref{thm:main}). 
Similar result for functions on $Sp(2,\mathbb{R})$ 
was proved by Hirano-Ishii-Oda~\cite{hio}, which motivated the study of this paper. 

\subsection*{Acknowledgement}
This paper is based on a lecture given by the author in ``Workshop on representation 
theory" held at 
Tambara Institute of Mathematical Sciences 
The University of Tokyo on October 8, 2006. 
The author thanks Professor Toshio Oshima and 
Professor Hisayoshi Matumoto 
for helpful discussions and comments. 

\section{Preliminaries}
\subsection{The Heckman and Opdam hypergeometric function}

In this subsection, we review on the Heckman-Opdam hypergeometric 
function associated with a root system. See \cite{heckman} 
and \cite{opdam} for details. 

Let $\mathfrak{a}$ be a Euclidean space of dimension $n$. 
For $\alpha\in \mathfrak{a}^*\setminus\{0\}$ define
\[
\alpha^{\vee}=\frac{2\alpha}{(\alpha,\alpha)}. 
\]
Let $R$ denote a reduced root system in $\mathfrak{a}^*$. 
Choose a positive system $R_+\subset R$ and 
let $B$ denote the set of the simple roots. 
Let $W$ denote 
the Weyl group for $R$. 
For $\alpha\in R$ let $k_\alpha$ be a non-negative number 
such that $k_{w\alpha}=k_\alpha$ for all $w\in W$. 
We call $k\,:\,\alpha\mapsto k_\alpha$ a multiplicity function. 
We put 
\[
\rho(k)=\frac12\sum_{\alpha\in R_+}k_\alpha \alpha.
\]
We often identify $\mathfrak{a}^*$ with $\mathfrak{a}$. 

Let $A=\exp \mathfrak{a}$ and 
\[
A_+=\{a\in A\,:\,\alpha(\log a)>0\text{ for all }\alpha\in R_+\}.
\]
Let $\{\xi_1,\dots,\xi_n\}$ be an orthonormal basis of $\mathfrak{a}$. 
Define
\begin{equation}
L(k)=\sum_{i=1}^n \partial_{\xi_i}^2+\sum_{\alpha\in R_+}
k_\alpha\frac{1+e^{-\alpha}}{1-e^{-\alpha}}\partial_\alpha.
\end{equation}
There exist a commutative algebra $\mathbb{D}(k)$ 
of differential operators containing $L(k)$ and an 
isomorphism $\gamma\,:\,\mathbb{D}(k)\rightarrow S(\mathfrak{a})^W$. 

Let $Q$ be the $\mathbb{Z}$-span of $R$ and 
$Q_+$ be the $\mathbb{Z}_+$-span of $R_+$. 
There exists a solution $\Phi(\lambda,k;a)$ for
\begin{equation}
\label{eqn:hgs}
D\varphi=\gamma(D)(\lambda)
\varphi
\end{equation}
of the form
\begin{equation}
\label{eqn:hcser}
\Phi(\lambda,k;a)=\sum_{\mu\in Q_+}
\Gamma_\mu(\lambda,k)e^{(\lambda-\rho(k)-\mu)(\log a)},
\quad \Gamma_0(\lambda,k)=0.
\end{equation}
The coefficients $\Gamma_\mu(\lambda,k)$ are determined 
by recurrence relations coming from $L(k)$. 

If $\lambda\in\mathfrak{a}_\mathbb{C}^*$ satisfies condition, 
\begin{equation}
\label{eqn:generic}
(2\lambda+\mu,\mu)\not=0\,\text{ for all }\, \mu\in Q\setminus\{0\}, 
\end{equation}
then $\{\Phi(w\lambda,k;a)\,:\,w\in W\}$ 
forms a basis of solution space of (\ref{eqn:hgs}) 
on $A_+$. 

Define
\begin{equation*}
\tilde{c}(\lambda,k)=\prod_{\alpha\in R_+}
\frac{\Gamma((\lambda,\alpha^{\vee}))}{\Gamma((\lambda,\alpha^{\vee})+k_\alpha)}
\end{equation*}
and
\begin{equation*}
c(\lambda,k)=\frac{\tilde{c}(\lambda,k)}{\tilde{c}(\rho(k),k)}.
\end{equation*}
Define
\begin{equation}
\label{eqn:hgfex}
F(\lambda,k;a)=\sum_{w\in W}{c}(w\lambda,k)\Phi(w\lambda,k;a).
\end{equation}
The function $F$ is called the Heckman-Opdam hypergeometric 
function for the root system $R$. It is well behaved 
compared to $\Phi$. 

\begin{Thm}[Heckman-Opdam]
$F(\lambda,k;a)$ is a unique $W$-invariant solution 
for (\ref{eqn:hgs}) that is  analytic in $a\in A$,  
holomorphic in $\lambda\in\mathfrak{a}_\mathbb{C}^*$, 
and 
\begin{align*}
& F(w\lambda,k;a)=F(\lambda,k;a) \quad (w\in W),\\
& F(\lambda,k;wa)=F(\lambda,k;a) \quad (w\in W).
\end{align*}
\end{Thm}

\subsection{Notation on Lie groups}

Let $G$ be a normal real form of a complex semisimple Lie group 
and $K$ a maximal compact subgroup. 
Let $\mathfrak{g}$ and $\mathfrak{k}$ be the Lie algebras 
of $G$ and $K$ respectively. Let $\theta$ denote the 
corresponding Cartan involution of $\mathfrak{g}$ and 
$\mathfrak{g}=\mathfrak{k}+\mathfrak{p}$ be the decomposition 
into $\pm 1$ eigenspaces of $\theta$. 
Equip the inner product $(\,\,,\,\,)$ on $\mathfrak{g}$ 
given by $(X,Y)=-B(X,\theta Y)$ $(X,\,Y\in\mathfrak{g})$ 
where $B(\,\,,\,\,)$ is the Killing form on $\mathfrak{g}$. 

Fix a maximal abelian subspace $\mathfrak{a}$ of $\mathfrak{p}$.  
Let $\Sigma=\Sigma(\mathfrak{g},\mathfrak{a})$ denote the set of 
the restricted roots.  Fix a positive system $\Sigma_+$ and 
let $\Pi$ denote the set of the simple roots. 
Notice that each root space has dimension 1, because 
we assume that $G$ is split. 
Put 
$\rho=\frac12\sum_{\alpha\in \Sigma_+}\alpha$. 
Let $W$ denote the Weyl group of $\Sigma$. 
It is isomorphic to $N_K(\mathfrak{a})/Z_K(\mathfrak{a})$, 
where $N_K(\mathfrak{a})$ (resp. $Z_K(\mathfrak{a})$)  
is the normalizer (resp. centralizer) of $\mathfrak{a}$ 
in $K$. 

Let $\mathfrak{n}$ be the sum of root spaces for the positive roots. 
Put 
$N=\exp \mathfrak{n}$ and $A=\exp \mathfrak{a}$. 
Then we have the Iwasawa decomposition $G=NAK=KAN$. 

\begin{Rem}
\label{rem:root}
In the previous section, we adopt notation of Heckman 
and Opdam. Relations between notation in the previous section and this 
section are given by
\begin{equation}
R=2\Sigma,\quad R_+=2\Sigma_+,\quad B=2\Pi,\quad 
 k_{2\beta}=\frac12m_{\beta}=\frac12\, \,\ (\beta\in \Sigma).
\end{equation}
 $L(k)$ is 
the radial part of the Laplace-Beltrami operator $L_{G/K}$ on 
$G/K$,  $\mathbb{D}(k)$ consists of 
the radial parts of invariant differential operators on 
$G/K$ with respect to the Cartan decomposition 
$G=KAK$, $c(\lambda,k)$ is the Harish-Chandra c-function,
 and the hypergeometric function $F(\lambda,k;a)$ is 
the radial part of the spherical function on $G/K$. 
\end{Rem} 

\subsection{The Whittaker function on a semisimple Lie group}
\label{subsec:toda}
In  this subsection, we review on the class one Whittaker functions on 
a split semisimple Lie group following Hashizume \cite{hashizume}. 

Let $\psi$ be a unitary character of $N$. 
We  denote the differential character of $\mathfrak{n}$ to 
$\sqrt{-1}\mathbb{R}$ by the same letter $\psi$. 
Let $C^\infty_\psi(G/K)$ denote the space of $C^\infty$-functions 
on $G$ satisfying $u(ngk)=\psi(n)u(g)$ for all 
$n\in N$, $g\in G$, and $k\in K$. By the Iwasawa decomposition 
the values of $u\in C^\infty_\psi(G/K)$ are completely 
determined by $u|_A$. 
Let $\mathbb{D}(G/K)$ denote the 
commutative algebra of left $G$-invariant differential 
operators on $G/K$ and $\chi_\lambda\,:\,
\mathbb{D}(G/K)\to \mathbb{C}$  the Harish-Chandra 
homomorphism. Let $\mathcal{A}_\psi(G/K,\mathcal{M}_\lambda)$ be the subspace 
of $C^\infty_\psi(G/K)$ defined by
\[
\mathcal{A}_\psi  (G/K,\mathcal{M}_\lambda) 
 =\{u\in C^\infty_\psi(G/K)\,:\, Du=\chi_\lambda(D)u\,\text{ for all }D\in
\mathbb{D}(G/K)\}. 
\]
Notice that $C^\infty_\psi(G/K,\mathcal{M}_\lambda)$ consists 
of real analytic functions, because $L_{G/K}$ is an elliptic 
differential operator. 

For $\beta\in\Pi$ let $X_\beta\in\mathfrak{g}^\beta$ be 
a unit root vector. For $\alpha\in B=2\Pi$ put 
$l_\alpha=-\sqrt{-1}\psi(X_{\alpha/2})$. 
For $u\in \mathcal{A}_\psi  (G/K,\mathcal{M}_\lambda) $, $\varphi=e^{-\rho}u|_A$ satisfies
\begin{equation}
\label{eqn:laptoda0}
\Big(\sum_{i=1}^n\partial_{\xi_i}^2
-2\sum_{\alpha\in B}l_\alpha^2 e^{\alpha}\Big)
\varphi=(\lambda,\lambda)\varphi.
\end{equation}
There exists a solution 
$\Psi_\text{T}(\lambda,\psi,a)$  for  
equation (\ref{eqn:laptoda0}) of the form
\begin{equation}
\label{eqn:psitoda}
\Psi_\text{T}(\lambda,\psi,a)=a^\lambda\sum_{\mu\in Q_+}
b_\mu(\lambda)a^\mu,\quad b_0(\lambda)=0.
\end{equation}
Moreover, extending function 
$u(a)=e^{\rho}\,\Psi_{\text{T}}(\lambda,\psi,a)$ on $A$ to $G$ 
so that $u\in C_\psi^\infty(G/K)$, 
it is also a joint eigenfunction of $\mathbb{D}(G/K)$ and 
belongs to  $\mathcal{A}_\psi(G/K,\mathcal{M}_\lambda)$. 
If  $\lambda\in\mathfrak{a}_\mathbb{C}^*$ satisfies 
condition (\ref{eqn:generic}), then 
$\{e^{\rho}\,\Psi_{\text{T}}(w\lambda,\psi,a)\,:\,w\in W\}$ 
forms a basis of $\mathcal{A}_\psi(G/K,\mathcal{M}_\lambda)|_A$ 
(cf. \cite[Corollary 5.3, Theorem 5.4]{hashizume}). 


For $\lambda\in\mathfrak{a}_{\mathbb{C}}^*$ define 
function $1_\lambda$ on $G$ by
\[
1_\lambda(nak)=a^{\lambda+\rho}\quad (n\in N,\,a\in A,\,k\in K).
\]
For $g\in G$ let $H(g)$ denote 
the element of $\mathfrak{a}$ defined by 
$g\in K\exp H(g)N$. We normalize the Haar measure 
$dn$ and $d\bar{n}$ on $N$ and $\bar{N}=\theta N$ by
\[
\theta(dn)=d\bar{n},\quad \int_{\bar{N}}e^{-2\rho(H(\bar{n}))}d\bar{n}=1,
\]
(cf. \cite[Ch. IV, \S 6]{helgason2}). 
Define
\begin{equation}
\label{eqn:jacquetint}
W(\lambda,\psi;g)=\int_N 1_\lambda(\bar{w}_0^{-1}ng)
\psi(n)^{-1}dn.
\end{equation}
Here $\bar{w}_0$ is a representative in $N_K(\mathfrak{a})$ of 
the longest element $w_0\in W$. 
For $\lambda\in\mathfrak{a}_\mathbb{C}^*$ with 
$\text{Re}\,\lambda>0$ (${}^\forall \,\alpha\in \Sigma_+$) 
the class one Jacquet integral $W(\lambda,\psi;g)$ converges 
absolutely and uniformly and belongs to $\mathcal{A}_\psi(G/K,\mathcal{M}_\lambda)$. 
Moreover $W(\lambda,\psi;g)$ can be continued 
to a meromorphic function of $\lambda\in\mathfrak{a}_\mathbb{C}^*$ 
as an element of $\mathcal{A}_\psi(G/K,\mathcal{M}_\lambda)$ 
 (cf. \cite[Theorem 6.6]{hashizume}). 
 We call this meromorphic continuation of 
$W(\lambda,\psi;a)$ the Whittaker function. 

The Whittaker function $W(\lambda,\psi;g)$ is up to a constant multiple a  
unique element of $\mathcal{A}_\psi(G/K,\mathcal{M}_\lambda)$ 
that is of moderate growth 
(cf. \cite[Theorem 9.1]{chm}). 

The analytic properties of the integral (\ref{eqn:jacquetint}) 
were studied by Jacquet~\cite{jacquet}, Schiffmann\cite{schiffmann}, 
Goodman-Wallach~\cite{gw1}, and Hashizume~\cite{hashizume0,  hashizume}, etc. 
$W(\lambda,\psi;g)$ satisfies the following functional equation.
\begin{equation}
\label{eqn:funceq}
W(\lambda,\psi;g)=M(w,\lambda,\psi)W(w\lambda,\psi;g)\quad (w\in W). 
\end{equation}
Here the function $M(w,\lambda,\psi)$ is given by the product 
formula
\begin{align}
\label{eqn:prodf1}
& M(ww',\lambda,\psi)=M(w',\lambda,\psi)M(w,w'\lambda,\psi)\quad 
(w,\,w'\in W), \\
\label{eqn:rankonef1}
& M(s_\alpha,\lambda,\psi)=
\left(\frac{2l_\alpha^2}{(\alpha,\alpha)} \right)^{(\lambda,\alpha^{\vee})}
\frac{\Gamma(-(\lambda,\alpha^\vee)+1/2)}{\Gamma((\lambda,\alpha^\vee)+1/2)}
\quad (\alpha\in B),
\end{align}
where $s_\alpha$ denote the simple reflection corresponding to $\alpha\in B$. 
(cf. \cite[(7.5)--(7.7)]{hashizume}. Notice again that 
$R=2\Sigma$ and $G$ is split.) 

Let $\mbox{\boldmath$c$}(\lambda)$ denote the Harish-Chandra $c$-function 
for the split Lie group $G$, which is given by 
$\mbox{\boldmath$c$}(\lambda)=c(\lambda,k)$ with 
$k_\alpha=1/2$ for all $\alpha\in R$. That is 
\begin{align*}
\mbox{\boldmath$c$}(\lambda) & 
=\frac{\tilde{\mbox{\boldmath$c$}}(\lambda)}{\tilde{\mbox{\boldmath$c$}}(\rho)},\\
{\tilde{\mbox{\boldmath$c$}}(\lambda)} & =
\prod_{\alpha\in R_+}\frac{\Gamma((\lambda,\alpha^\vee))}
{\Gamma((\lambda,\alpha^\vee)+\frac12)}.
\end{align*}
Hashizume~\cite[Theorem 7.8]{hashizume}  
expressed the Whittaker function $W(\lambda,\psi;a)$ as a 
linear combination of $\Psi_{\text{T}}(w\lambda,\psi,a)$ explicitly. 

\begin{Thm}[Hashizume]
\label{thm:hashizume}
Let $\psi$ be a non-degenerate character of $N$ and assume that  
$\lambda\in \mathfrak{a}_\mathbb{C}^*$ satisfies (\ref{eqn:generic}). 
Then 
\begin{equation}
\label{eqn:hashizumemain1}
W(\lambda,\psi;a)=a^\rho\sum_{w\in W}M(w_0w,\lambda,\psi)
\mbox{\boldmath$c$}({w}_0w\lambda)\Psi_{\text{\rm T}}(w\lambda,\psi;a)
\quad (a\in A_+).
\end{equation}
\end{Thm}

\section{Limit transition  from 
the Heckman-Opdam hypergeometric function to 
the Whittaker function}

\subsection{Limit transition from the Calogero-Moser Hamiltonian to the Toda Hamiltonian}
\label{subsection:limiteq}
In this subsection, we review on the limit transition from 
the quantum trigonometric Calogero-Moser model to the Toda model. 

Define a function $\delta(k)=\delta(k;a)$ by
\[
\delta(k)^{1/2}=\prod_{\alpha\in R_+}(e^{\frac12\alpha}-e^{-\frac12\alpha})^{k_\alpha}.
\]
We have 
\begin{align}
\label{eqn:sutherland}
\delta(k)^{1/2}\circ\{L(k) & +(\rho(k),\rho(k))\}\circ \delta(k)^{-1/2} \notag \\
& =\sum_{i=1}^n \partial_{\xi_i}^2+
\sum_{\alpha\in R_+}\frac{k_\alpha(1-k_\alpha)(\alpha,\alpha)}{4\sinh^2\frac12\alpha}
\end{align}
We denote the right hand side of (\ref{eqn:sutherland}) by $H_\text{CM}(k)$. 
It is the Hamiltonian for the trigonometric Calogero-Moser model.

Recall that $R=2\Sigma$ (Remark~\ref{rem:root}) and 
let  $B=2\Pi$ be the simple system of $R_+=2\Sigma_+$. 
We assume that $l_\alpha=1$ for all $\alpha\in B$. That is 
we assume that $\psi$ is a special non-degenerate unitary character of $N$. 
The left hand side of  (\ref{eqn:laptoda0}) gives 
the Hamiltonian for the quantum Toda model
\begin{equation}
\label{eqn:htoda2}
H_{\text{T}}
=\sum_{i=1}^n \partial_{\xi_i}^2-2
\sum_{\alpha\in B}
e^{\alpha}.
\end{equation}

Let $M$ be a positive real number. 
Define a positive multiplicity function 
$k_M$  by 
\begin{equation*}
k_M(\alpha)(k_M(\alpha)-1)(\alpha,\alpha)=2
e^{2M}
\end{equation*}
and define $a_M\in A$ by
\begin{equation*}
\log a_M=w_0\log a+M
\rho^\vee,
\end{equation*}
where $w_0$ is the longest element of $W$. 
Notice that 
\begin{equation*}
\rho^\vee=\frac12\sum_{\beta\in\Sigma_+}\beta^\vee
=\sum_{\alpha\in R_+}\alpha^\vee
\end{equation*}
is the Weyl vector of $\Sigma^\vee=2 R^\vee$ 
and $(\alpha,\rho^\vee)=1$ for all $\alpha\in\Pi=\frac12 B$ 
(cf. Bourbaki \cite[Ch~VI Proposition 29]{bourbaki}). 

We shall consider limits of the hypergeometric function when $M\to\infty$. 
Taking a limit of $H_{\text{CM}}(k)$, 
we have the following lemma. 
\begin{Lem}
\label{lem:limit1}
For any  $\varphi\in C^\infty(A)$, 
\begin{equation}
\lim_{M\to\infty}H_{CM}(k_M)\,\varphi(a_M)=H_{\text{\rm T}}
\,\varphi(a).
\end{equation}
\end{Lem}

This limit procedure was proved by Inozemtsev \cite{inozemtsev} 
(see also \cite[Section 7]{etingof} and \cite{takasakietal}). 

\subsection{Limit transition of eigenfunctions}

Define 
\begin{equation*}
\Psi_\text{CM}(\lambda,k;a)=\delta(k;a)^{1/2}\,\Phi(\lambda,k;a).
\end{equation*}
By (\ref{eqn:hcser}) and (\ref{eqn:sutherland}), 
$\varphi(a)=\Psi_\text{CM}(\lambda,k;a)$ is 
of the form
\begin{equation}
\label{eqn:hhg}
\Psi_\text{CM}(\lambda,k;a)=\sum_{\mu\in \Lambda}b_\mu(\lambda,k)
e^{(\lambda-\mu)(\log a)}, \quad b_0(\lambda,k)=1
\end{equation}
and it is a solution of 
\begin{equation}
\label{eqn:cme}
H_\text{CM}(k)\,\varphi=(\lambda,\lambda)\,\varphi. 
\end{equation}

On the other hand, as we have seen in subsection~\ref{subsec:toda}, 
there is  
a series solution $\varphi(a)=\Psi_{\text{T}}(\lambda;a)$ of
\begin{equation*}
H_{\text{T}}
\,\varphi=(\lambda,\lambda)\,\varphi.
\end{equation*}

\begin{Prop}
\label{prop:limit2}
If $\lambda\in\mathfrak{a}_\mathbb{C}^*$ satisfies condition (\ref{eqn:generic}), then 
\begin{equation}
\label{eqn:limhc}
\lim_{M\to\infty}e^{-(\lambda,\rho^\vee)M}
\Psi_{\text{\rm CM}}(\lambda,k;a_M)
=\Psi_{\text{\rm T}}(w_0\lambda;a)\quad (a\in A_+).
\end{equation}
The convergence is uniform on each subchamber
\[
\{a\in A_+\,:\,\alpha(\log a)>c>0\,\,(\alpha\in B)\},
\]
where $c>0$ is arbitrary. 
\end{Prop}
\begin{proof}
The proof is an easy modification of the estimate of the Harish-Chandra 
series due to Gangolli~\cite{gangolli} (see also Helgason~\cite[Ch IV \S 5]{helgason2}). 

Substituting $k_M$ and $a_M$ 
to  (\ref{eqn:hhg}) and (\ref{eqn:cme}), 
we have
\begin{equation}
\label{eqn:seriesM}
\Psi_{\text{\rm CM}}(\lambda,k;a_M)=e^{(\lambda,\rho^\vee)M}
\sum_{\mu\in Q_+}\tilde{b}_\mu(\lambda,M)e^{(w_0\lambda+\mu)(\log a)}, 
\quad \tilde{b}_0(\lambda,M)=1
\end{equation}
and 
\begin{align}
\label{eqn:deM}
\Big(\sum_{i=1}^n\partial_{\xi_i}^2-2\sum_{\alpha\in R_+}
e^{2M}
\sum_{j=1}^\infty j\,  e^{-j((\alpha,\rho^\vee)M+w_0\alpha(\log a))} & 
\Big)
\Psi_\text{CM}(\lambda,k_M;a_M) \notag \\
& =(\lambda,\lambda)\Psi_\text{CM}(\lambda,k_M;a_M)
\end{align}
Substituting (\ref{eqn:seriesM}) to (\ref{eqn:deM})
 gives the recurrence relation
\begin{align}
 (2w_0\lambda+ & \mu,\mu)\tilde{b}_\mu(\lambda,M) \notag \\
& = 2\sum_{\alpha\in R_+}
\sum_{j\geq 1,\,\mu+jw_0\alpha\in Q_+}
e^{(2-j(\alpha,\rho^\vee))M}\,j\,\tilde{b}_{\mu+jw_0\alpha}(\lambda,M).
\label{eqn:rec1}
\end{align}
Since $(\alpha,\rho^\vee)=2$ for $\alpha\in B=2\Pi$ and 
$(\alpha,\rho^\vee)\geq 4$ for $\alpha\in R_+\setminus B$, 
 recurrence relation (\ref{eqn:rec1}) converges to 
\begin{equation}
\label{eqn:rec2}
 (2w_0\lambda+  \mu,\mu)\tilde{b}_\mu(\lambda,\infty)
=2\sum_{\alpha\in B}
\tilde{b}_{\mu-\alpha}(\lambda,\infty)
\end{equation}
as $M\to\infty$. (\ref{eqn:rec2}) is nothing but the 
recurrence relation for the coefficients in  
expansion (\ref{eqn:psitoda}) for $\Psi_\text{T}(w_0\lambda;a)$ 
(cf. \cite[\S 4]{hashizume}). 

For $\mu\in Q_+$ we write $\mu=\sum_{\alpha\in B} n_\alpha\alpha$ and 
put $n(\mu)=\sum_{\alpha\in B} n_\alpha$. 
Choose a constant $c$ such that 
\[
|2(w_0\lambda+\mu,\mu)|\geq c\, n(\mu)
\]
for all $\mu\in Q_+$. 
By (\ref{eqn:rec1}) we have
\begin{align*}
|\tilde{b}_\mu(\lambda,M)|  
& \leq 2c^{-1}\sum_{\alpha\in R_+}
\sum_{j\geq 1,\,\mu+jw_0\alpha\in Q_+}
e^{(2-(j\alpha,\rho^\vee))M}\,j\,|\tilde{b}_{\mu+jw_0\alpha}(\lambda,M)| \\
& \leq 2c^{-1}\sum_{\alpha\in R_+}
\sum_{j\geq 1,\,\mu+jw_0\alpha\in Q_+}
j\,|\tilde{b}_{\mu+jw_0\alpha}(\lambda,M)|
\end{align*}
for $M>0$. We can prove in the same 
way as the proof of \cite[Ch IV Lemma 5.3, Lemma 5.6]{helgason2} that 
there exists a constant $K_{a,l}$ such that
\[
|\tilde{b}_\mu(\lambda,M)|\leq K_{\lambda,a}a^\mu
\]
for all $\mu\in\Lambda$. 
This estimate shows the convergence of the 
series (\ref{eqn:seriesM}) and also guarantees the 
limit transition (\ref{eqn:limhc}). 
\end{proof}

Now we state and prove our main result:

\begin{Thm}
\label{thm:main}
Assume that 
$\lambda\in \mathfrak{a}^*_\mathbb{C}$ satisfies (\ref{eqn:generic}). 
Then 
\begin{align}
\label{eqn:main}
\lim_{M\to\infty} & 
\delta(k;a_M)^{1/2}\,\tilde{c}(\rho(k),k)
\prod_{\alpha\in R_+}
\Gamma(k_M(\alpha))\,F(\lambda,k_M;a_M) \notag \\
& =\tilde{\mbox{\boldmath$c$}}(\rho)f(\lambda)\,
a^{-\rho}\,W(\lambda,\psi;a), 
\end{align}
where
\begin{equation}
f(\lambda)=\prod_{\alpha\in R_+}
\left(\frac{(\alpha,\alpha)}{2
}\right)^{(\lambda,\alpha^\vee)/2}
\Gamma((\lambda,\alpha^\vee)+\tfrac12)
\end{equation}
and $\psi$ is a unitary character of $N$ defined by $l_\alpha=1\,
(\alpha\in B)$. 
\end{Thm}

\begin{proof}
By the following formula for the Gamma function
\[
\lim_{x\to\infty}\frac{\Gamma(\mu+x)}{\Gamma(x)\,x^\mu}=1,
\]
we have 
\begin{equation}
\label{eqn:limitcf}
\tilde{c}(\lambda,k_M)\sim 
f(\lambda)
\tilde{\mbox{\boldmath$c$}}(\lambda)
\prod_{\alpha\in R_+}
\frac{e^{-(\lambda,\alpha^\vee)M}}{\Gamma(k_M(\alpha))}
\end{equation}
as $M\to\infty$.
By (\ref{eqn:hgfex}), Proposition~\ref{prop:limit2}, 
and (\ref{eqn:limitcf}), we have
\begin{align}
\lim_{M\to\infty} & 
\delta(k;a_M)^{1/2}\tilde{c}(\rho(k),k)
\prod_{\alpha\in R_+}
\Gamma(k_M(\alpha))\,F(\lambda,k_M;a_M) \notag \\
& =\sum_{w\in W}
f(w\lambda)\tilde{\mbox{\boldmath$c$}}(w\lambda)
\Psi_\text{T}(w_0w\lambda;a).
\label{eqn:asympt1}
\end{align}

On the other hand, by Theorem~\ref{thm:hashizume}, 
the right hand side of (\ref{eqn:main}) is a linear combination 
of $\Psi_\text{T}(w\lambda;a)\,\,(w\in W)$, where 
the coefficient of $\Psi_\text{T}(w\lambda;a)$ is given by
\begin{equation}
d(w,\lambda):=
f(\lambda)
M(w_0w,\lambda,\psi)\tilde{\mbox{\boldmath$c$}}(w_0w\lambda).
\end{equation}
For $\beta\in R$, it follows from (\ref{eqn:prodf1}) and (\ref{eqn:rankonef1}) that 
\begin{align*}
d  (w,s_\beta\lambda)  & =f(s_\beta\lambda)
M(w_0w,s_\beta\lambda,\psi)\tilde{\mbox{\boldmath$c$}}(w_0ws_\beta\lambda) \\
& = 
f(s_\beta\lambda)
M(s_\beta,\lambda,\psi)^{-1}M(w_0ws_\beta,\lambda,\psi)
\tilde{\mbox{\boldmath$c$}}(w_0ws_\beta\lambda) \\
&= d(ws_\beta,\lambda).
\end{align*}
The last equality follows from $s_\beta(B\setminus\{\beta\})=B\setminus\{\beta\}$. 
Thus the right hand side of  (\ref{eqn:main})  is $W$-invariant with 
respect to $\lambda$. We have 
$d({w_0},\lambda)  =
f(\lambda)
\tilde{\mbox{\boldmath$c$}}(\lambda)$ and it coincides 
with the coefficient of $\Psi_\text{T}(\lambda,l;a)$ in the 
right hand side of (\ref{eqn:asympt1}), hence the result follows.  
\end{proof}

\begin{Ex}
For $R$ of type $A_1$
\begin{equation*}
F(\lambda,k;a_t)={}_2 F_1(\tfrac12(k-\lambda),
\tfrac12(k+\lambda)\,;\,k+\tfrac12\,;\,-\sinh^2 t), 
\end{equation*}
where ${}_2 F_1(a,b,c;z)$ is the Gauss hypergeometric function. 
Theorem~\ref{thm:main} reads
\begin{equation*}
\lim_{k\to\infty}k^{-1/2}2^{-k}\sinh^k(-t+M)
F(\lambda,k;a_{-t+M})=\frac{1}{\sqrt{\pi}}K_\lambda(e^t/2),
\end{equation*}
where $K_\lambda(z)$ is the Macdonald function. 
\end{Ex}

\begin{Rem}
We restrict ourselves to split semisimple Lie groups, 
because the Hamiltonian (\ref{eqn:htoda2}) of the Toda lattice 
depends only on reducible root system. 
Class one Whittaker function given by the Jacquet integral 
for a non-split semisimple Lie group is a constant multiple 
of the Whittaker function for the split Lie group of the 
same indivisible restricted roots. 

We can change 
parameters $l_\alpha^2\,(\alpha\in B)$ in the left hand side of 
(\ref{eqn:laptoda0}) by making 
a shift of variables, as it was pointed out by \cite[\S 2.1]{takasakietal}. 
Let $\varpi_\alpha\,(\alpha\in B)$ 
denote the fundamental weights corresponding to $B$. If we put 
\[
\log a=\log a'-\sum_{\alpha\in B,\ l_{\alpha}\not=0}
\frac{2}{(\alpha,\alpha)}\varpi_\alpha \log l_{\alpha}^2, 
\]
then 
\[
\alpha(\log a)=\left\{
\begin{matrix}
\alpha(\log a')-\log l_{\alpha}^2 \quad (\alpha\in B,\,
l_\alpha\not=0) \\
\alpha(\log a') \quad (\alpha\in B,\,l_\alpha=0)
\end{matrix}
\right.
\]
Thus in the new coordinates 
$\sum_{i=1}^n\partial_{\xi_i}^2-
2\sum_{\alpha\in B} l_\alpha^2\,e^{\alpha}$ 
becomes
$
\sum_{i=1}^n\partial_{\xi_i}^2-
2\sum_{\alpha\in B,\ l_\alpha\not=0} e^{\alpha}
$. 
Moreover, if $l_\alpha=0$ for some $\alpha$, then the Whittaker function 
can be reduced to lower rank cases. 
This is the reason why we assume $l_\alpha=1$ for all $\alpha\in B$. 
\end{Rem}

\address{Department of Applied Mathematics \\
Okayama University of Science \\
Okayama 700-0005
 \\Japan}
{shimeno@xmath.ous.ac.jp}
\end{document}